\documentclass[12pt]{amsart}
\usepackage{amsmath,amsthm,amscd,color}
\usepackage{amssymb}
\usepackage{amsfonts}
\usepackage{latexsym}
\usepackage{mathrsfs}
\usepackage{hyperref}
\usepackage{geometry}
\usepackage{fancyvrb}
\usepackage{dsfont}
\usepackage{pgfplots}
\usepackage[normalem]{ulem}
\usepgfplotslibrary{fillbetween}

\newcommand{\eps}{\varepsilon}
\newcommand{\nequiv}{\not\equiv}
\renewcommand{\(}{\left(}
\renewcommand{\)}{\right)}
\allowdisplaybreaks

\newtheorem*{theoremA}{Theorem 1}
\newtheorem*{theoremB}{Theorem 2}
\newtheorem*{theoremC}{Theorem 3}
\newtheorem{theorem}[equation]{Theorem}

\newtheorem{lemma}[equation]{Lemma}

\newtheorem{proposition}[equation]{Proposition}
\newtheorem*{corollary}{Corollary}

\newtheorem{conjthm}[equation]{Conjecture/Theorem}
\theoremstyle{remark}
\newtheorem{definition}[equation]{Definition}

\newtheorem*{remark}{\bf Remark}

\numberwithin{equation}{section}
\numberwithin{table}{section}

\DeclareMathOperator{\re}{Re}

\DeclareMathOperator{\SL}{SL}

\renewcommand{\(}{\left(}
\renewcommand{\)}{\right)}
\newcommand{\ts}{\textstyle}

\newcommand{\Z}{{\mathbb{Z}}}
\newcommand{\R}{{\mathbb{R}}}
\newcommand{\C}{{\mathbb{C}}}

\newcommand{\uH}{{\mathbb{H}}}

\newcommand{\smatr}[4]{\left( \begin{smallmatrix} #1 & #2 \\ #3 & #4 \end{smallmatrix}\right)}

\newcommand{\abcd}{\smatr{a}{b}{c}{d}}

\author[Jiang, Rolen, Woodbury]{Yuze Jiang, Larry Rolen, and Michael Woodbury}
\address{Department of Mathematics\\Columbia University\\ 2990 Broadway \\ New York, NY 10027, USA}
\email{yj2595@columbia.edu} 
\address{Department of Mathematics \\
Vanderbilt University \\
1420 Stevenson Center \\
Nashville, TN 37240, USA}
\email{larry.rolen@vanderbilt.edu}
\address{Department of Mathematics \\Brown University \\ 151 Thayer St \\ Providence, RI 02912, USA}
\email{michael\_woodbury@brown.edu}
\title{Generalized Frobenius Partitions, Motzkin paths, and Jacobi Forms}
\date{\today}

\begin{document}
\maketitle

\begin{abstract}
We show how Andrews' generating functions for generalized Frobenius partitions can be understood within the theory of Eichler and Zagier as specific coefficients of certain Jacobi forms.  This reformulation leads to a recursive process which yields explicit formulas for the generalized Frobenius partition generating functions in terms of infinite $q$-products.   In particular, we show that specific examples of our result easily reestablish previously known formulas, and we describe new congruences, both conjectural and proven, in additional cases.  The modular structure of Jacobi forms indicates that \emph{all} of the coefficients of the forms are of interest.  We give a combinatorial definition of these ``companion series'' and explore their combinatorics via the counting of Motzkin paths.
\end{abstract}

\section*{Introduction}

There is a long and storied history of studying the partition function $p(n)$, which counts the number of \emph{partitions} of $n$, i.e., nonincreasing sequences of positive integers summing to $n$.  A particularly pivotal development in the theory of partitions was the Hardy-Ramanujan-Rademacher Circle Method which gives the exact formula 
\begin{equation}\label{eq:HRR}
 p(n)=\frac{1}{\pi\sqrt{2}}\sum\limits_{k=1}^{\infty}A_k(n)\sqrt{k}\frac{d}{dn}\left(\big(n-1/24\big)^{-1/2}\sinh\left(\frac{2\pi}{k\sqrt{6}}\big(n-1/24\big)^{1/2}\right)\right),
\end{equation}
where $A_k(n)$ is a certain sum of roots of unity.

A key observation in the Hardy--Ramanujan--Rademacher Circle Method is that the generating function for the partition numbers
 \[ C\Phi_1(q):=\sum_{n=0}^\infty p(n)q^n = \prod_{j=1}^\infty \frac{1}{1-q^j} =: \frac{1}{(q;q)_\infty} = \frac{q^{\frac{1}{24}}}{\eta(\tau)} \qquad \mbox{(where $q:=e^{2\pi i\tau}$)} \]
is very nearly a \emph{modular form}.  Indeed, the Dedekind eta function $\eta(\tau)$ is a modular form with multiplier of weight $\frac12$.  The proof of \eqref{eq:HRR} relies heavily on the fact that $\eta(\tau)$ is modular.  Although $p(n)$ has a purely combinatorial description/definition, by couching the partition numbers as coefficients of a modular form, one can introduce very powerful analytic techniques to study $p(n)$.  

In this paper, we will apply the principle of making use of the analytic (i.e., modular) properties of certain generating functions to glean information about their coefficients.  In particular, we will show that the coefficients of a certain class of \emph{Jacobi forms} $F_k(z;\tau)$ have an interpretation as generating functions of combinatorially interesting objects.  These coefficients encompass and generalize $C\Phi_1(q)$ and the well-known generating functions $C\Phi_k(q)$ of $k$-colored Frobenius partitions.  We make use of the Jacobi form structure to reveal a hidden recursive relation among the various generating functions. Our results serve as a sample of applying standard techniques from the theory of Jacobi forms to combinatorial generating functions, including coefficients of certain infinite products, such as those which arise from applications of Andrews' General Principle described below.  Since such formulas appear throughout partition theory and combinatorics, it would be interesting to search for additional examples like these and work out similar results.  Our results also show that these coefficients can be described in terms of Motzkin paths, certain lattice paths.

In order to define the functions $F_k(z;\tau)$, we begin with the \emph{Jacobi theta function}:
\begin{equation}\label{eq:theta}
 \vartheta(z;\tau) := \sum_{n\in \frac{1}{2}+\Z}e^{\pi i n^2\tau + 2\pi i n\left(z+\frac{1}{2}\right)},
\end{equation}
where $z,\tau\in \C$ with $\re(\tau)>0$.  Then, for each $k\in \Z_{> 0}$, we set
\begin{equation}\nonumber
F_k(z;\tau) := \left(\frac{-\vartheta\left(z+\frac12;\tau\right)}{q^{\frac1{12}}\eta(\tau)}\right)^k.
\end{equation}
See Section~\ref{sec:JacobiForms} for background on Jacobi forms.  In particular, note that a Jacobi form $F(z;\tau)$ must be periodic in both variables: $F(z+1;\tau)=F(z;\tau+1)=F(z;\tau)$, and is holomorphic in both $\tau$ and $z$.  Hence there is an expansion in terms of $q=e^{2\pi i \tau}$ and $\zeta:=e^{2\pi i z}$.

In Section~\ref{sec:combinatorics} we establish the following.
\begin{theoremA}\label{thm:mainthmA}
The $\zeta^a$ coefficient of $F_k(z;\tau)$ is 
\begin{equation}\label{eq:CPsikaDef}
 C\Psi_{k,a}(q):= \sum_{n=0}^\infty c\psi_{k,a}(n)q^n,
\end{equation}
where $c\psi_{k,a}(n)$ is the number of $(k,a)$-colored Frobenius partitions of weight $n$ (see Definition~\ref{def:kaFrobparts} or, equivalently, Definition~\ref{def:kaFrob-Drake}).  In particular, $C\Psi_{k,k/2}(q)=C\Phi_k(q)$.
\end{theoremA}

\begin{remark}
Andrews introduced the notation
 \[ C\Phi_k(q) := \sum_{n=0}^\infty c\phi_k(n) q^n \]
where $c\phi_k(n)$ is the number of $k$-colored Frobenius partitions of $n$ (see Section~\ref{sec:combinatorics}).
\end{remark}
\begin{remark}
It is interesting to note that the Jacobi forms whose Fourier coefficients we study here are related to a well-studied important family. Specifically, up to a change of variables and a power of $\eta(\tau)$, the $\zeta$-Fourier coefficients of the Jacobi form $\vartheta(z+1/2;\tau)^m/\vartheta(z;\tau)^n$ with $m\geq2$ and $n\geq1$ are essentially important affine Lie superalgebra characters known as Kac-Wakimoto characters \cite{KW,KW94} (see Section 20.2 of \cite{HMFBook}). These have been studied in a number of recent papers (see, e.g., \cite{BCR14, BF14, BFM15, BO09, F11}), where they gave rise to objects such as mock modular and quantum modular forms. Given this history, it is interesting that here we add a study of the modular properties and applications of the analogous family where $n=0$ and the Jacobi form in question is holomorphic.
\end{remark}

The functions $F_k(z;\tau)$ were first addressed by the last two authors and Kathrin Bringmann in \cite{BRW16}, in which the fact that the $\zeta^{\frac{k}{2}}$ coefficient of $F_k(z;\tau)$ is $C\Phi_k(q)$ is utilized.  With Kathrin Bringmann's generous permission, the current paper is an update of that work.  Unbeknownst to us at that time, in \cite{Drake09} Brian Drake considered the \emph{constant} term, which he denoted $C\Psi_k(q)$ (in the case that $k=2\ell$ is even).  Many of our methods in Sections~\ref{sec:combinatorics} and \ref{sec:motzkin} came about as a result of trying to refashion the work in \cite{Drake09} to deal with all of the coefficients.  A key difference is that although Definition~\ref{def:kaFrob-Drake} (a direct generalization of Drake's definition) is advantageous for describing the link to Motzkin paths, our equivalent formulation Definition~\ref{def:kaFrobparts} is much simpler to state.  Moreover, Theorem~1 follows as an immediate consequence of Definition~\ref{def:kaFrobparts} and a subtle modification to George Andrews' so-called ``General Principle'' from \cite{Andrews84} (see Lemma~\ref{lem:ExtendedGP} below).  The equivalence of these definitions and applications to Motzkin paths are discussed in Section~\ref{sec:motzkin}.

Given the fact the the partition function satisfies so many striking congruence relations, it is natural to ask whether $c\phi_k(n)$ (or $c\psi_{k,a}(n)$) also exhibit simple congruences.  Indeed, the answer is yes, and there is a long history of results of this type.  In \cite{Andrews84} it was shown that $c\phi_p(n)\equiv 0 \pmod{p^2}$ for primes $p\nmid n$.  These general results where extended in \cite{K2} and further in \cite{Sellers}.  The recent paper \cite{GS2014} proved that there exist several infinite families of congruences, for example $c\phi_{5N+1}(5n+4)\equiv 0\pmod{5}$ for any $n,N\in \Z_{>0}$. 

Using explicit realizations of $C\Phi_k$ in terms of $q$-series for small values of $k$, many authors have found additional congruences.  The case $k=3$ was studied in \cite{Kolitsch}, $k=4$ in \cite{BS11,Lin,Sellers4,Xia,ZWW17}, and $k=6$ in \cite{BS15,H}. The techniques of these papers are quite similar in that they rely on having an explicit $q$-series representation for $C\Phi_k$ in terms of the $q$-Pochhammer symbol.  The idea employed by Andrews to establish these $q$-series formulas for $C\Phi_2$ and $C\Phi_3$ (See \eqref{eq:Andrews2} and \eqref{eq:Andrews3}) have been largely mirrored on a case-by-case basis.  This same general strategy was also used, for example, by Drake to establish 
 \[ C\Psi_{2,0}(q) = \frac{2}{(q;q)_\infty (q;q^4)_\infty(q^2;q^4)_\infty^2(q^3;q^4)_\infty}, \]
where $(a;q)_\infty := \prod\limits_{n=1}^\infty \left(1-aq^{n-1}\right)$ is the usual $q$-Pochhammer symbol.  This procedure, however, becomes increasingly tedious as $k$ grows.  In this paper we show that by using the language of Jacobi forms, there is a simple recursive formula which gives a formula for any specific $k$. Using our recursive formula, the previously known formulas for the generating functions $C\Phi_k(q)$ and $C\Psi_{k,0}(q)$ can be recovered.  We demonstrate this for several cases in Section~\ref{sec:examples}.  We remark that in \cite{CWY} a method for finding $C\Phi_k(q)$ that uses modular forms on a case-by-case basis is presented.  

As described in Section~\ref{sec:JacobiForms} any Jacobi form such as $F_k(z;\tau)$ has a theta decomposition. In the case that $k=2\ell$ is even, we have
 \[ F_{2\ell}(z;\tau) = \displaystyle{\sum_{a\pmod{2\ell}}H_{\ell,a}(\tau)\vartheta_{\ell,a}(z;\tau)}, \]
where
\begin{equation}\label{eq:Jthetaam}
 \vartheta_{m,a}(z;\tau) := \sum_{\substack{r\in \Z\\ r\equiv a\pmod{2m}}} q^{\frac{r^2}{4m}}\zeta^r = \sum_{n\in \Z} q^{\frac{(2mn+a)^2}{4m}} \zeta^{2mn+a}.
\end{equation}
Our next main result gives an iterative formula for the functions $H_{\ell,c}(\tau)$ in terms of $$\theta_{m,a}(\tau):=\vartheta_{m,a}(0;\tau)=\sum_{n\in \Z}q^{\frac{(2mn+a)^2}{4m}}.$$  This can then be used to calculate $C\Psi_{k,a}(q)$ for any $k$ and $a$. We remark that considering all of the coefficients is essential to making the recursion work.  Note that this approach differs from Andrews' treatment which only looked at a single term.

\begin{theoremB}
For any $j\in \Z$, let
\begin{equation}\label{eq:base}
  h_{1,j}(\tau):=\theta_{1,j+1}(\tau).
\end{equation}
Given $\{h_{\ell,c}(\tau): 0\leq c < 2\ell\}$, recursively define for each $0\leq b\leq \ell$
\begin{align}\label{eq:inducttwo}
 \begin{split}
 h_{\ell+1,b}(\tau) & := h_{1,b}(\tau)\theta_{\ell(\ell+1),b\ell}(\tau)h_{\ell,0}(\tau)
       +h_{1,b-\ell}(\tau)\theta_{\ell(\ell+1),b\ell-\ell(\ell+1)}(\tau)h_{\ell,\ell}(\tau) \\
  & \quad + \sum_{c=1}^{\ell-1} h_{1,b-c}(\tau)
    \big( \theta_{\ell(\ell+1),b\ell-c(\ell+1)}(\tau) 
          +\theta_{\ell(\ell+1),b\ell+c(\ell+1)}(\tau) \big) h_{\ell,c}(\tau),
 \end{split}
\\ \label{eq:inductone}
  h_{\ell+\frac12,b+\frac12}(\tau) & := \sum_{c\pmod{2\ell}} h_{\ell,c}(\tau)\theta_{\ell(2\ell+1),c(2\ell+1)-\ell(2b+1)}(\tau).
\end{align}
For $\ell < b \leq 2\ell$, set $h_{\ell,b}(\tau):=h_{\ell,2\ell-b}(\tau)$.  Then, for all $\ell\in \Z_{>0}$ and $\delta\in\{0,1\}$, 
\begin{equation}\label{equation2ldelta}
 (-1)^{\delta}\vartheta\Big(z+\frac12;\tau\Big)^{2\ell+\delta} = \sum_{b \pmod{2\ell+\delta}} h_{\ell+\frac{\delta}{2},b+\frac{\delta}{2}}(\tau)\vartheta_{\ell+\frac{\delta}{2},b+\frac{\delta}{2}}(z;\tau). 
\end{equation}
\end{theoremB}

\begin{remark}
The recursive formula that gives $h_{\ell+\frac12,b}(\tau)$ in terms of  $h_{\ell,b}(\tau)$ holds even if $\ell\in \frac12\Z$, but one must sum over half integers.  By further making the substitions $\ell\mapsto \ell+1/2$, $b\mapsto b-1/2$, and $c\mapsto c+1/2$, \eqref{eq:inductone} gives the formula
\begin{equation}\nonumber
   h_{\ell+1,b}(\tau) = \sum_{c\pmod{2\ell+1}} h_{\ell+\frac12,c+\frac12}(\tau)\theta_{(2\ell+1)(\ell+1),(2c+1)(\ell+1)-(2\ell+1)b}(\tau)
\end{equation}
which holds for $\ell\in\Z_{>0}$.
\end{remark}

\begin{corollary}\label{thetacor}
All Fourier coefficients of $\vartheta(z+1/2;\tau)^{k}$ can be given as combinations of Dedekind eta functions and Klein forms (see \eqref{eq:Klein}). 
In particular, there exists an algorithm to compute $C\Psi_{k,a}(q)$ as a sum of products of $q$-Pochhammer symbols.
\end{corollary}
\begin{remark}
As we shall see in the proof of this corollary in Section~\ref{sec:proofofCor}, the ``central'' coefficient $h_{\ell+\frac\delta2,\ell+\frac\delta2}$ in the theta decomposition \eqref{equation2ldelta} is, up to a simple multiple of a power of $q$ and a power of $\eta$, equal to  $C\Phi_k(\tau)$.
\end{remark}

The proof of Theorem~2 and its corollary are the subject of Section~\ref{sec:proofs}.

In addition to giving a robust method for deriving additional formulas, we expect that considering $C\Phi_k$ from the point of view of Jacobi forms provides additional means by which congruences can be studied.  As an example of how ``modularity'' has been used previously, see \cite{Ono} (later extended in \cite{L}), which established congruence properties for $\phi_3(n) \pmod{7}$ by relating $\sum_{n\geq 1}\phi_3(9n)q^n$ to modular forms.  Based on numerical observation, we collect several congruence properties of $c\psi_{k,a}(n)$ in Section~\ref{sec:congruences}, some of which we have proved.

Our final result realizes the coefficients of the Jacobi forms $F_k(z;\tau)$, hence $(k,a)$-colored F-partitions, in terms of Motzkin paths, which are certain lattice paths.  Our work extends that of Drake in \cite{Drake09} in which a connection between Motzkin paths and the constant term of $F_{2\ell}(z;\tau)$ was first observed.  In order to state our result, we set the following notation.

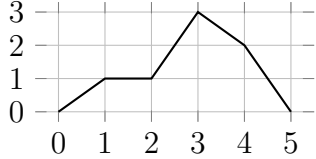
\begin{figure}[h]
\begin{center} \begin{tikzpicture}
  \begin{axis}[
  width=0.35\textwidth,
       height=0.21\textwidth,
xtick={0,1,...,5},
ytick={0,1,...,3},
    separate axis lines,
    y axis line style= { draw opacity=0 },
    x axis line style= { draw opacity=0 },
    ymajorgrids,
		 xmajorgrids,
	   ]		
    \addplot[name path=E,mark=none,thick,sharp plot] coordinates 
		{(0,0) (1,1) (2,1) (3,3) (4,2) (5,0)};
  \end{axis}
\end{tikzpicture} \end{center}
\caption{A $0$-shifted Motzkin path of length $5$ and rank $2$}\label{fig:MP}
\end{figure}

\begin{definition}\label{def:MotzkinPath}
For $k\in \Z_{> 0}$ and $a\in \frac{k}{2}+\Z$ ($a\geq 0$), 
an \emph{$a$-shifted Motzkin path of length $n$ and rank $\frac{k}{2}$} is a path from $(0,0)$ to $(n,a(n-1))$ which does not go below the $x$-axis and consists of steps $(1,j)$ with $j\in \{a-\frac{k}{2},a-\frac{k}{2}+1,\ldots, a+\frac{k}{2}\}$.  A \emph{$(k,a)$-Motzkin path} is an $a$-shifted Motzkin path of rank $\frac{k}{2}$ for which the $j$-th step is colored by a subset $T\subseteq \{1,2,\cdots, k\}$ of size $\# T= \frac{k}{2}-a+j$.
\end{definition}
The $0$-shifted Motzkin paths correspond to the usual definition of a Motzkin path.  In the case that $a=0$ and $k=2\ell$ is even, one obtains what Drake refers to as a \emph{colored Motzkin path with length $n$ and rank $\ell$} (See \cite[p.~3948]{Drake09}).  An example of a $(4,0)$-Motzkin path of length $5$ (without coloring) is given in Figure~\ref{fig:MP}. 

Given an $a$-shifted Motzkin path $w$, let $a_M(w)$ denote the area between the Motzkin path and the $x$-axis, and we set
 \[ \widetilde{CM}_{k,a}^n(q):=\sum q^{a_M(w)}, \]
where the sum is over all $(k,a)$-Motzkin paths $w$ of length $n$.  Let
 \[ CM_{k,a}^{2n+1}(q):=q^{\frac{k}{2}n(n+1)+an(2n+1)}\widetilde{CM}_{k,a}^{2n+1}\big(q^{-1}\big). \]
\begin{theoremC}
With the definitions given above, 
\begin{equation*}
    \lim_{n\to \infty}CM_{k,a}^{2n+1}(q)=C\Psi_{k,a}(q).
\end{equation*}
\end{theoremC}
In Section~\ref{sec:motzkin} we prove that Definition~\ref{def:kaFrobparts} and  Definition~\ref{def:kaFrob-Drake} of $(k,a)$-colored $F$-partitions are equivalent, and then we proceed to prove Theorem~3. 

\section*{Acknowledgments}

The authors are extremely grateful to Kathrin Bringmann who allowed us to take the preprint \cite{BRW16} as a starting point for this work.  We also thank Chung Hang (Kevin) Kwan and Ben Goldman who met regularly with the first and third authors via Zoom throughout the Summer of 2020.  Kevin was instrumental in helping the first author to learn about $q$-series, and Ben patiently checked the work in Section~\ref{sec:congruences}. Furthermore, we are grateful to Bruce Berndt for helpful conversations in particular pointing out Schroeter's formula to the authors, we thank Jeffrey Lagarias for useful comments which improved the exposition of the paper, and we thank Ian Wagner for useful conversations on the topics of the paper.  The second author is grateful for support from a grant from the Simons Foundation (853830, LR), support from a Dean’s Faculty Fellowship from Vanderbilt University, and to the Max Planck Institute for Mathematics in Bonn for its hospitality and financial support.

\section{Combinatorics}\label{sec:combinatorics}

Given a partition $n$, one obtains its associated Frobenius coordinates by reading off the leg and arm lengths along the diagonal of its Ferrers diagram (see \cite{Andrews84} for complete definitions of these terms). This gives a bijection between partitions of $n$ and arrays
\begin{align}\label{eq:array}
 \left( \begin{array}{cccc}
    a_1 & a_2 & \cdots & a_m \\
    b_1 & b_2 & \cdots & b_m \end{array} \right)
\end{align}
with integral coordinates $a_1>a_2>\text{\dots}>a_m\geq 0$ and $b_1>b_2>\text{\dots}>b_m\geq 0$ such that $n=m+\sum_{j=1}^m(a_j+b_j)$. For example, associated to the partition $\lambda=(5,4,4,2,1)$ of $16$, we have the following diagram, with the boxes along the diagonal labeled with the symbol $\bullet$: 

\begin{figure}[h]
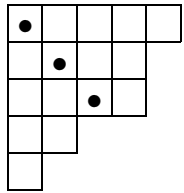
\begin{center}
{\def\lr#1{\multicolumn{1}{|@{\hspace{.6ex}}c@{\hspace{.6ex}}|}{\raisebox{-.3ex}{$#1$}}}
\raisebox{-.6ex}{$\begin{array}[b]{*{5}c}\cline{1-5}
\lr{\bullet}&\lr{\phantom{x}}&\lr{\phantom{x}}&\lr{\phantom{x}}&\lr{\phantom{x}}\\\cline{1-5}
\lr{\phantom{x}}&\lr{\bullet}&\lr{\phantom{x}}&\lr{\phantom{x}}\\\cline{1-4}
\lr{\phantom{x}}&\lr{\phantom{x}}&\lr{\bullet}&\lr{\phantom{x}}\\\cline{1-4}
\lr{\phantom{x}}&\lr{\phantom{x}}\\\cline{1-2}
\lr{\phantom{x}}\\\cline{1-1}
\end{array}$}
}\end{center}
\caption{The Ferrers diagram for the partition $\lambda=(5,4,4,2,1)$}
\end{figure}

The leg length for any diagonal cell is then the number of boxes below it, and the arm length is the number of cells to its right. So the Frobenius coordinates for $\lambda$ are 
\begin{align*}
 \left( \begin{array}{ccc}
    4 & 2 & 1 \\
    4 & 2 & 0\end{array} \right)
    .
\end{align*}

A \emph{generalized Frobenius partition} or \emph{F-partition} is an array as in \eqref{eq:array} but where the rows are allowed to come from more general sets.  In \cite{Andrews84}, Andrews introduced this notion and considered two particular types of F-partitions.  In this paper we are particularly concerned with one of these types, namely \emph{the generalized Frobenius partitions in $k$-colors}.  This is obtained by requiring that the sequences $(a_j)_{1\leq j\leq m}$ and $(b_j)_{1\leq j\leq m}$ as in \eqref{eq:array} are each strictly decreasing sequences of integers selected from $k$ copies of the nonnegative integers.  By \emph{strictly decreasing}, we mean with respect to the following lexicographical ordering: if $m$ belongs to the $j$th copy of $\Z_{\geq0}$, which we denote by writing $m=m_j$, and $n=n_\ell$ belongs to the $\ell$th copy of $\Z_{\geq 0}$, then $m_j<n_\ell$ precisely when $m<n$ or $m=n$ and $j<\ell$.

Following Andrews' notation, we denote the number of such partitions of $n$ by $c\phi_k(n)$ and define the generating function
\begin{equation}\label{eq:CPhik}
 C\Phi_k(q) := \sum_{n=0}^\infty c\phi_k(n)q^n.
\end{equation}
The special case of $k=1$ corresponds to the usual partition function. 

Andrews computed $C\Phi_2(q)$ and $C\Phi_3(q)$ using his ``General Principle'' which we describe and generalize below.  In particular, he proved that
\begin{equation}\label{eq:Andrews2}
 C\Phi_2(q)=\frac{\left(q^2;q^4\right)_\infty}{\left(q;q^2\right)_\infty^4\left(q^4;q^4\right)_\infty}
\end{equation}
and
\begin{equation}\label{eq:Andrews3}
 C\Phi_3(q)
= \frac{\left(q^{12};q^{12}\right)_\infty\left(q^{6};q^{12}\right)_\infty^3}{\left(q;q^6\right)_\infty^5\left(q^5;q^6\right)_\infty^5\left(q^4;q^4\right)_\infty^2\left(q^3;q^6\right)_\infty^7} + 4q\frac{\left(q^{12};q^{12}\right)_\infty\left(q^4;q^4\right)_\infty}{\left(q^6;q^{12}\right)_\infty\left(q^2;q^4\right)_\infty\left(q;q\right)_\infty^3}.
\end{equation}

We generalize the notion of $k$-colored F-partitions as follows.
\begin{definition}\label{def:kaFrobparts}
Given $k\in \Z_{> 0}$ and $a\in \Z+\frac{k}{2}$ ($a\geq 0$), a \emph{$(k,a)$-colored F-partition} is a two rowed array
\begin{equation} \label{eq:kaFrobPartitionnew}
\begin{pmatrix}
a_1 & a_2 & \cdots  & a_r \\
b_1 & b_2 & \cdots  & b_s \\
\end{pmatrix}
\end{equation}
satisfying the following conditions.
\begin{itemize}
 \item[(i)] Each entry $a_\ell,b_j$ belongs to one of $k$ copies of nonnegative integers.
 \item[(ii)] Each row is decreasing with respect to the lexicographical ordering.
 \item[(iii)] The pair $(r,s)\neq(0,0)$ of non-negative integers satisfies
\begin{equation}\label{eq:sconditionsnew}
r-s=a-\frac{k}{2}.
\end{equation}
\end{itemize}

The \emph{weight} of a $(k,a)$-colored F-partition as in \eqref{eq:kaFrobPartitionnew} is
\begin{equation}
n:=r+\sum\limits_{\ell=0}^{r}a_\ell+\sum\limits_{j=0}^{s}b_j.
\end{equation}
We let $c\psi_{k,a}(n)$ be the number of $(k,a)$-Frobenius partitions of weight $n$, and we set $C\Psi_{k,a}(q):=\sum\limits_{n=0}^\infty c\psi_{k,a}(n)q^n$ to be its generating function.
\end{definition}

Note that if $a=\frac{k}{2}$, then by \eqref{eq:sconditionsnew} $r=s$, and Definition~\ref{def:kaFrobparts} recovers the $k$-colored Frobenius partitions.  In other words, $C\Psi_{k,\frac{k}{2}}(q) = C\Phi_k(q)$.

Our main result with regards to $(k,a)$-colored F-partitions is Theorem~1, namely that they comprise the coefficients of $F_k(z;\tau)$.  Following Andrews \cite{Andrews84}, we have the following key result.

\begin{lemma}\label{lem:ExtendedGP}
Let $P_A(m,n)$ be the number of ordinary partitions of $n$ into $m$ parts subject to a restriction $A$, and let $f_A(z)=\sum P_A(m,n)z^mq^n$ be the associated generating function.  Then the generating function $C\Phi_{A,B}(q)$ for the two-rowed arrays \eqref{eq:kaFrobPartitionnew} whose first row satisfies restriction $A$ and second row satisfies $B$ with difference of length equal to $a-\frac{k}{2}$ is the constant term of 
\begin{align*}
z^{\frac{k}{2}-a}f_A(zq)f_B(z^{-1}) 
\end{align*}
\end{lemma}

The proof relies on the following result, which, when specialized to the case $a=\frac{k}{2}$, is precisely Andrews' ``General Principle'' referenced in the introduction.

\begin{lemma}\label{lem:CPsiInterpretnew}
The function $C\Psi_{k,a}(q)$ is the constant term (with respect to $\zeta$) of 
\begin{equation}
\zeta^{\frac{k}{2}-a}\prod\limits_{n=0}^{\infty}(1+\zeta q^{n+1})^k(1+\zeta^{-1}q^{n})^k, 
\end{equation}
with $a \in  \mathbb{Z}+\frac{k}{2}$.  Moreover, if we define
 \[ F_k(z;\tau) := \left(\frac{-\vartheta(z+\frac12;\tau)}{q^{\frac{1}{12}}\eta(\tau)} \right)^k. \]
then $C\Psi_{k,a}(q)$ is also the $\zeta^a$ coefficient of $F_k$.
\end{lemma}

\begin{proof}
We seek to prove that $C\Psi_{k,a}(q)$ is the $\zeta^a$ coefficient  of $F_k(z;\tau)$
\begin{align*}
    F_k(z;\tau)=&\left(\frac{-\vartheta(z+\frac{1}{2};\tau)}{q^{\frac{1}{12}}\eta(\tau)}\right)^k
    =\zeta^\frac{k}{2}\prod\limits_{n=0}^{\infty}(1+\zeta q^{n+1})^k(1+\zeta^{-1}q^n)^k.
\end{align*}
The latter equality here is immediate by the Jacobi triple product formula \eqref{JTP}.  The above, in turn, is the constant term of
\begin{equation*}
\zeta^{\frac{k}{2}-a}\prod\limits_{n=0}^{\infty}(1+\zeta q^{n+1})^k(1+\zeta^{-1}q^n)^k,
\end{equation*}
which is precisely the generating function for $C\Psi_{k,a}(q)$ as in Lemma~\ref{lem:ExtendedGP}.
\end{proof}

\section{Background on Jacobi forms}\label{sec:JacobiForms}

\subsection{Jacobi forms}\label{JacobiDefnSection}
We first require the definition of Jacobi forms, whose theory was laid out by Eichler--Zagier \cite{EZ}.  A \emph{holomorphic Jacobi form of weight $\kappa$ and index $m$ on $\SL_2(\Z)$} is a holomorphic function $\varphi\colon\C\times \uH\to \C$ ($\uH$ is the complex upper half-plane) which satisfies the following conditions. 
\begin{itemize}
 \item[(i)] For all $\gamma=\abcd\in \SL_2(\Z)$, 
 \[ \varphi\left(\frac{z}{c\tau+d};\frac{a\tau+b}{c\tau+d}\right)=(c\tau+d)^\kappa e^{\frac{2\pi i mcz^2}{c\tau+d}}\varphi(z;\tau). \]
 \item[(ii)] For all $\lambda,\mu\in \Z$,
 \[ \varphi(z+\lambda\tau+\mu;\tau)=e^{-2\pi i m\left(\lambda^2\tau+2\lambda z\right)}\varphi(z;\tau). \]
 \item[(iii)]  The function $\varphi$ has a Fourier expansion of the form
 \[ \varphi(z;\tau)= \sum_{n,r\in \Z} c(n,r)q^n\zeta^r\]
with $c(n,r)=0$ unless $4mn\geq r^2$.
\end{itemize}
\begin{remark}
As with the ordinary theory of modular forms, there are suitable modifications of this definition for congruence subgroups of $\SL_2(\Z)$, half-integral weight (and index), as well as multiplier systems (as arise in examples such as Jacobi's theta function below). For ease of exposition, we omit these technical definitions here, opting to only give the necessary transformations for the basic Jacobi form needed here, namely the Jacobi theta function.
\end{remark}

\subsection{The Jacobi theta function}\label{sec:JacTheta}
The Jacobi form which is of most importance in this paper is the Jacobi theta function defined in \eqref{eq:theta}.  It satisfies the well-known transformation properties (see, e.g., \cite{HMFBook})
\begin{equation*}
\vartheta(z;\tau+1) = e^{\frac{\pi i}{4}} \vartheta(z;\tau),\qquad
\vartheta\left(-\textstyle{\frac{z}{\tau};-\frac{1}{\tau}}\right) = i\sqrt{-i\tau}e^{\frac{\pi i z^2}{\tau}}\vartheta(z;\tau),
\end{equation*}
and is an example of a holomorphic Jacobi form of weight $1/2$ and index $1/2$.

Moreover, the Jacobi theta function satisfies the well-known Jacobi triple product identity 
\begin{align}\label{JTP}
 \vartheta(z;\tau) =  -iq^{\frac18}\zeta^{-\frac12}\left(q;q\right)_\infty\left(\zeta q^{-1};q\right)_\infty\left(\zeta^{-1};q\right)_\infty.
\end{align}

\subsection{Theta decomposition}\label{sec:thetadecomp} 

The main structural result on Jacobi forms that we exploit in this paper is that every holomorphic Jacobi form can be expressed in terms of the Jacobi theta functions $\vartheta_{m,a}(z;\tau)$, which we recall were given in \eqref{eq:Jthetaam} as 
\begin{equation*}
 \vartheta_{m,a}(z;\tau) = \sum_{\substack{r\in \Z\\ r\equiv a\pmod{2m}}} q^{\frac{r^2}{4m}}\zeta^r = \sum_{n\in \Z} q^{\frac{(2mn+a)^2}{4m}} \zeta^{2mn+a}.
\end{equation*}
In the case that $m$ and $a$ are half-integers, which is relevant for this paper, we take the latter sum as the definition.

The following theorem shows that every Jacobi form has a theta decomposition.  

\begin{theorem}[Eichler-Zagier] \label{th:EZ}
Suppose $\varphi(z;\tau)$ is a holomorphic Jacobi form of weight $\kappa$ and index $m$.  Then 
 \[ \varphi(z;\tau) = \sum_{a\pmod{2m}}h_a(\tau)\vartheta_{m,a}(z;\tau), \]
where $(h_a(\tau))_{a\ ({\rm mod}\ 2m)}$ is a vector-valued modular form of weight $\kappa-1/2$.
\end{theorem}

So, in particular, a modified version of Theorem~\ref{th:EZ} covering a slightly more general class of Jacobi functions to which $F_k(z;\tau)$ belongs implies that $F_k(z;\tau)$ has a similar decomposition; Theorem~2 provides this decomposition directly.  Note that we include Theorem~\ref{th:EZ} as a motivating principle, but not as essential to our proof.  Indeed, Theorem~2 is deduced directly and independently of Theorem~\ref{th:EZ}.

For future reference, note that the modular forms $\theta_{m,a}(\tau)=\vartheta_{m,a}(0;\tau)$ satisfy the relations
\begin{equation}\label{eq:thetamarelation}
 \theta_{m,2mk\pm a}(\tau)= \theta_{m,a}(\tau),
\qquad \mbox{and}\qquad
 \theta_{m,a}(k\tau)= \theta_{km,ka}(\tau),
\end{equation}
for all $k\in \Z_{>0}$.

\subsection{Klein forms}\label{sec:Klein}
Further distinguished modular forms are obtained by specializing the Jacobi theta function to torsion points of $\C/(\Z+\tau\Z)$ and dividing by a power of the $\eta$-function.  To be more precise, for $a\in \mathbb{R}$ we have the \emph{Klein forms} 
\begin{equation}\label{eq:Klein}
t_{a,0}(\tau):=-q^{\frac{a^2}{2}-\frac{a}{2}+\frac{1}{12}}\frac{\left(q^a;q\right)_\infty\left(q^{1-a};q\right)_\infty}{\left( q;q\right)_\infty^2}.
\end{equation}
Specifically, for $a\in \mathbb{Q}$, these functions are modular forms of weight $-1$, holomorphic in the upper half-plane and with possible poles and zeros only at cusps.  (See \cite{KL}.)

\section{Proof of Theorem~2 and its Corollary}\label{sec:proofs}

\subsection{First examples of Theorem~2}\label{sec:firstexamples}
This section is devoted to the cases $k=1,2$, a detailed analysis of which sheds light on the general procedure. The Fourier expansion
\begin{equation}\label{Theta12}
-\vartheta\(z+\frac12;\tau\) = \sum_{n\in \Z} q^{\frac12\(n+\frac12\)^2}\zeta^{n+\frac12}
\end{equation}
is an immediate consequence of \eqref{eq:theta}.  Therefore, 
\begin{equation}\nonumber
 F_1(z;\tau) = \sum_{n\in \Z}\left( \frac{q^{\frac12\(n+\frac12\)^2}}{q^{\frac{1}{12}}\eta(\tau)} \right) \zeta^{n+\frac12},
\end{equation}

We can rederive Andrews' formula for $C\Phi_2(q)$ similarly by employing the following result of \cite{BM2013}, where the authors used the notation $A(z;\tau):=\vartheta(z;\tau)^2/\eta(\tau)^6$.  Accounting for the shift $z\mapsto z+1/2$, they gave the following lemma. 

\begin{lemma}\label{lem:theta2}
The square of $\vartheta(z+\frac12;\tau)$ has the theta decomposition
 \[ \vartheta\left(z+\textstyle{\frac12};\tau\right)^2 = \theta_{1,1}(\tau)\vartheta_{1,0}(z;\tau)+\theta_{1,0}(\tau)\vartheta_{1,1}(z;\tau).\]
\end{lemma}

\begin{proof}
The following proof follows that in \cite{BM2013} but for the reader's convenience, we give it here.  First, using \eqref{Theta12}, note that
 \[  \vartheta\left(z+\textstyle{\frac12};\tau\right)^2 = \sum_{r,s\in\Z}q^{\frac12\(\(r+\frac12\)^2+\(s+\frac12\)^2\)}\zeta^{r+s+1}. \]
Thus, making the change of variables $R:=r+s$ and $S:=r-s$, we obtain
\begin{equation*}
 \vartheta\left(z+\textstyle{\frac12};\tau\right)^2 = \sum_{\substack{R,S\in \Z\\ R\equiv S\pmod{2}}}q^{\frac{1}{4}\left((R+1)^2+S^2\right)}\zeta^{R+1}.
\end{equation*}
Depending on the parity of $R$ and $S$, this sum splits naturally into two parts.  The contribution from $R$ and $S$ even equals 
 \[ \sum_{S\in \Z} q^{S^2}\sum_{R\in \Z}q^{\left(R+\frac12\right)^2}\zeta^{2R+1} = \theta_{1,0}(\tau)\vartheta_{1,1}(z;\tau). \]
The contribution from $R$ and $S$ odd is 
 \[ \sum_{S\in \Z} q^{\left(S+\frac12\right)^2} \sum_{R\in \Z}q^{R^2}\zeta^{2R} = \theta_{1,1}(\tau)\vartheta_{1,0}(z;\tau). \]
Adding these together, the result follows.
\end{proof}

Two well-known consequences of the Jacobi triple product formula are
 \begin{align}\label{eq:theta0eta}
   \theta_{1,0}(\tau) =\ & \frac{\eta(2\tau)^5}{\eta(\tau)^2\eta(4\tau)^2} = \frac{\left(q^2;q^{2}\right)_\infty^5}{\left(q;q\right)_\infty^2\left(q^4;q^{4}\right)_\infty^2},\\
\nonumber
   \theta_{1,1}(\tau) =\ & 2\frac{\eta(4\tau)^2}{\eta(2\tau)} = 2q^{\frac14} \frac{\left(q^{4};q^{4}\right)_\infty^2}{\left(q^{2};q^{2}\right)_\infty}. 
 \end{align}
 In order to find the coefficient of $\zeta$ in $F_2(z;\tau)$, we use \eqref{eq:theta0eta} and Lemma~\ref{lem:theta2}.  This leads to
\begin{align*}
 C\Phi_2(q) = \frac{q^{\frac14}\theta_{1,0}(\tau)}{q^{\frac16}\eta(\tau)^2}\ & = \frac{\left(q^2;q^{2}\right)_\infty^5}{\left(q;q\right)_\infty^4\left(q^4;q^{4}\right)_\infty^2},
\end{align*}
which agrees with \eqref{eq:Andrews2}.

\subsection{Preliminary results}

We next prove two lemmas, both of which can be viewed as variations on Lemma~\ref{lem:theta2}.  These lemmas describe how one can obtain the theta decomposition of the weight and index $1$ Jacobi form (with multiplier) 
\begin{equation}\label{eq:fortwomorelem}
 \vartheta_{1,\eps}(z;\tau)\vartheta_{\ell,c}(z;\tau) 
\end{equation}
for $\eps\in\{0,1\}$, and for 
\begin{equation}\nonumber
 \vartheta\left( z+\frac12;\tau \right) \vartheta_{\ell,c}(z;\tau),
\end{equation}
respectively.  

With Lemma~\ref{lem:theta2} in hand, knowing the theta decomposition of \eqref{eq:fortwomorelem} with $\ell=1$ and $c\in\{0,1\}$ one can find the theta decomposition of $\vartheta(z+\frac12;\tau)^3$.  More generally, for any $\ell\in\Z_{>0}$, Lemma~\ref{lem:theta2} provides the key step to go from the theta decomposition of $\vartheta(z+\frac12;\tau)^\ell$ to that of $\vartheta(z+\frac12;\tau)^{\ell+2}$.

\begin{lemma}\label{lem:onemore}
Assume the notation as above.  Then, for $c,\ell\in \Z_{\geq0}$ with $\ell>c $, 
 \[ -\vartheta\left( z+\frac12;\tau \right) \vartheta_{\ell,c}(z;\tau) 
    = \sum_{a\pmod{2\ell+1}} \theta_{\ell(2\ell+1),c-2\ell a-\ell}(\tau)\vartheta_{\ell+\frac12,a+c+\frac12}(z;\tau). \]
\end{lemma}

\begin{proof}
As in the proof of Lemma~\ref{lem:theta2}, we see that
\begin{align}\nonumber
 -\vartheta\left( z+\frac12;\tau \right) \vartheta_{\ell,c}(z;\tau) 
  = & \sum_{r,s\in \Z} q^{ \frac12\left( s+\frac12\right)^2+\frac{(2\ell r+c)^2}{4\ell}} \zeta^{2\ell r+s+c+\frac12}.
\end{align}
We make the change of variables $r= \frac{R+S}{2\ell+1}, s = \frac{S-2\ell R}{2\ell+1}$. As $(r,s)$ runs through $\Z^2$, $(R,S)$ run through those elements in $\Z^2$ satisfying, as claimed, $R+S\equiv 0\pmod{2\ell+1}$.  Thus we find 
\begin{align*}
 -\vartheta&\left( z+\frac12;\tau \right) \vartheta_{\ell,c}(z;\tau) \\
  = & \sum_{\substack{R,S\in \Z \\ R+S\equiv 0 \pmod{2\ell+1}}} q^{ \frac{\ell}{2\ell+1}\left( R+\frac{c}{2\ell}-\frac12 \right)^2 + \frac{1}{2(2\ell+1)}\left( S+c+\frac12 \right)^2 }\zeta^{S+c+\frac12} \\
  = & \hskip-6pt \sum_{a\ \mathrm{mod}\, 2\ell+1} \hskip-2pt \left(
 \sum_{R\in \Z} q^{\frac{\ell}{2\ell+1}\left((2\ell+1)R-a+\frac{c}{2\ell}-\frac12\right)^2} 
\sum_{S\in\Z} q^{\frac{1}{2(2\ell+1)}\left( (2\ell+1)S+a+c+\frac12 \right)^2} \zeta^{(2\ell+1)S+a+c+\frac12} \right).
\end{align*}
From this, we conclude the claim.
\end{proof}

\begin{lemma}\label{lem:theta1eps}
With the same conditions as in Lemma \ref{lem:onemore}, we have that
 \[ \vartheta_{1,\eps}(z;\tau)\vartheta_{\ell,c}(z;\tau) = \sum_{a \pmod{\ell+1}}\theta_{\ell(\ell+1),(2a+\eps)\ell-c}(\tau)\vartheta_{\ell+1,2a+c+\eps}(z;\tau). \]
\end{lemma}

The proof of Lemma~\ref{lem:theta1eps} follows the same idea as that of Lemma \ref{lem:onemore}.  We leave it as an exercise for the interested reader.

\begin{remark}
After proving Lemma~\ref{lem:theta1eps}, Bruce Berndt informed the authors of a very similar result, known as Schroeter's formula, which gives the identity
 \[ T\left(x,q^a\right)T\left(x,q^b\right) = \sum_{k=0}^{a+b-1}y^kq^{bk^2}T\left(xyq^{2bk},q^{a+b}\right)T\left(y^ax^{-b}q^{2abk},q^{ab(a+b)}\right), \]
where $T(x,q):=\sum_{n\in \Z}x^n q^{n^2}$. Since making an earlier version of the paper \cite{BRW16} first available on the arXiv, the authors have since learned that the main theorem of \cite{Cao} is a generalization of Schroeter's theorem and that Lemmas~\ref{lem:theta2}--\ref{lem:theta1eps} should follow from it as well.
\end{remark}

\subsection{Proof of Theorem~2}

Let $\ell\in \Z_{>0}$ and $\delta\in \{0,1\}$.  Our method of proof is to give an iterative procedure for obtaining $\vartheta(z+\frac12;\tau)^{2\ell+1+\delta}$ from $\vartheta(z+\frac12;\tau)^{2\ell}$.  Suppose that we have a decomposition
\begin{equation}\label{equation2l}
 \vartheta\Big(z+\frac12;\tau\Big)^{2\ell} = \sum_{c \pmod{2\ell}} h_{\ell,c}(\tau)\vartheta_{\ell,c}(z;\tau)
\end{equation}
with $h_{\ell,c}=h_{\ell,2\ell-c}$.  By Lemma~\ref{lem:theta2}, \eqref{equation2l} is true in the case $\ell=1$.  This establishes \eqref{eq:base}.  We must then show that defining $h_{\ell+\frac{1-\delta}{2},b+\frac{\delta}{2}}$ as in \eqref{eq:inducttwo} and \eqref{eq:inductone}, it is also true that \eqref{equation2ldelta} holds.

We first treat the case $\delta=1$.  In this case, applying Lemmas~\ref{lem:theta2} and \ref{lem:theta1eps}, we find that 
\begin{align*}
 \vartheta &\Big(z+\frac12;\tau\Big)^{2(\ell+1)}
   =  \Big(\theta_{1,1}(\tau) \vartheta_{1,0}(z;\tau)+\theta_{1,0}(\tau)\vartheta_{1,1}(z;\tau)\Big) \vartheta\Big(z+\frac12;\tau\Big)^{2\ell} \\
    &=  \sum_{c\pmod{2\ell}}h_{\ell,c}(\tau)
    \left( \theta_{1,1}(\tau)\sum_{a\pmod{\ell+1}} \theta_{\ell(\ell+1),2a\ell-c}(\tau)\vartheta_{\ell+1,2a+c}(z;\tau)\right. \\
   &\qquad \qquad \qquad \qquad \quad
  + \left. \theta_{1,0}(\tau)\sum_{\alpha\pmod{\ell+1}} \theta_{\ell(\ell+1),(2\alpha+1)\ell-c}(\tau)\vartheta_{\ell+1,2\alpha+1+c}(z;\tau)\right).
\end{align*}

For a given $b\pmod{2\ell+2}$, we now collect all of the terms in which $\vartheta_{\ell+1,b}(z;\tau)$ appears.  In the sum over $a$, there is exactly one such term if $c\equiv b\pmod{2}$ and none otherwise.  This term comes from the unique $a$ such that $2a+c\equiv b \pmod{2\ell+2}$.  Similarly, in the sum over $\alpha$, we get one term exactly if $c\nequiv b\pmod{2}$, namely for the unique $\alpha$ for which $2\alpha+c+1\equiv b\pmod{2\ell+2}$.  Hence,
\begin{align*}
 h_{\ell+1,b}(\tau) = &
\sum_{\substack{c\pmod{2\ell}\\ c\equiv b \pmod{2}}} 
  \theta_{1,1}(\tau) \theta_{\ell(\ell+1),(b-c)\ell-c}(\tau)h_{\ell,c}(\tau) \\
 & \quad +
\sum_{\substack{c\pmod{2\ell}\\ c\nequiv b \pmod{2}}}
  \theta_{1,0}(\tau) \theta_{\ell(\ell+1),(b-c)\ell-c}(\tau)h_{\ell,c}(\tau). 
\end{align*}
Since $h_{1,b-c}$ is equal to $\theta_{1,1}$ or $\theta_{1,0}$ exactly depending on whether $b$ and $c$ have the same parity or not, this can be simplified to
\begin{equation}\label{eq:inductalmost}
 h_{\ell+1,b}(\tau) = \sum_{c\pmod{2\ell}}
  h_{1,b-c}(\tau) \theta_{\ell(\ell+1),b\ell-c(\ell+1)}(\tau)h_{\ell,c}(\tau).
\end{equation}
From this formula and the inductive hypothesis that $h_{\ell,b}=h_{\ell,2\ell-b}$, it follows that $h_{\ell+1,b}=h_{\ell+1,2\ell+2-b}$.  This and the fact that 
$\theta_{\ell(\ell+1),a}$ depends only on $a\hskip -1pt\pmod{2\ell(\ell+1)}$ implies that \eqref{eq:inductalmost} now simplifies readily to \eqref{eq:inducttwo}.

The case of $\delta=0$ corresponds to \eqref{eq:inductone}.  This case follows the same logic as above with the only difference being to employ Lemma~\ref{lem:onemore} instead of Lemma~\ref{lem:theta1eps}.  We leave the details to the reader.

\subsection{Proof of Corollary to Theorem~2}\label{sec:proofofCor}

The corollary to Theorem~2 is a direct consequence of the fact that the function $\theta_{m,b}$ can be expressed as an infinite product in terms of $q$-Pochhammer symbols. Indeed, from \eqref{eq:theta}, \eqref{JTP}, and \eqref{eq:Klein}, it follows that
\begin{align}
 \theta_{m,b}(\tau) = & q^{\frac{b^2}{4m}}\left(q^{2m};q^{2m}\right)_\infty\left(q^{m-b};q^{2m}\right)_\infty\left(q^{m+b};q^{2m}\right)_\infty \nonumber \\
 = & -q^{\frac{m}{12}}\left( q^{2m};q^{2m}\right)_\infty^3 t_{\frac{1}{2}+\frac{b}{2m},0}(2m\tau). \label{eq:thetambasprod}
\end{align}

We next claim that, up to a multiple by a power of $q$ and a power of $\eta$, $C\Phi_k(q)$ is $h_{\frac{k}{2},\frac{k}{2}}$. First, note that 
 
\[ C\Phi_k(q)=C\Psi_{k,\frac k2}(q)=[\zeta^{\frac k2}]F_{k}(z;\tau)=(-1)^kq^{-\frac{k}{12}}\eta(\tau)^{-k}[\zeta^{\frac k2}]\vartheta(z+\frac12;\tau)^k, \]
where $[\zeta^a]f(z;\tau)$ denotes the $\zeta^a$ coefficient of $f(z;\tau)$. By \eqref{equation2ldelta}, the coefficients in the theta decomposition of $\vartheta(z+\frac12;\tau)^k$ are (up to a sign) the $h_{k,r}$'s. A standard fact in the theory of Jacobi forms (see for example Section~3.1 of \cite{ChengDuncan}) is that for a Jacobi form $f(z;\tau)$ of index $m$, the coefficients in the theta decomposition are equal to the Fourier coefficients in $\zeta$ up to $q$ powers. Specifically, in the language of Theorem~\ref{th:EZ}, we have
\[
q^{-\frac{r^2}{4m}}[\zeta^r]f(z;\tau)=h_r(\tau).
\]
Combining these observations shows that $C\Phi_k$ is indeed a multiple of a power of $q$ and a power of $\eta$ times $h_{\frac k2,\frac k2}$.

Since Theorem~2 expresses  $h_{\frac{k}{2},\frac{k}{2}}$ as a combination of the functions $\theta_{m,b}$, the corollary is evident.

\section{Special cases}\label{sec:examples}

In this section, whenever we are only dealing with functions of the single variable $\tau$, or equivalently $q$, in order to give cleaner formulas we drop the arguments of the function if possible except in the statements of the theorems.

\subsection{The case of $C\Phi_6$}\label{sec:cphi6}
The identities of \eqref{eq:thetamarelation} are used repeatedly in the following formulas.

Applying \eqref{eq:inducttwo} of Theorem~2 twice, we find first that 
\begin{equation}\label{eq:h2b}
  h_{2,0} = \theta_{1,1}^2\theta_{2,0}+\theta_{1,0}^2\theta_{2,2}, \quad
  h_{2,1} = 2\theta_{1,0}\theta_{1,1}\theta_{2,1}, \quad
  h_{2,2} = \theta_{1,1}^2\theta_{2,2}+\theta_{1,0}^2\theta_{2,0},
\end{equation}
and then that 
\begin{align}
  h_{3,0}  
  = &\ \theta_{1,1}\theta_{6,0}h_{2,0}+2\theta_{1,0}\theta_{6,3}h_{2,1}
   + \theta_{1,1}\theta_{6,6}h_{2,2}, 
\\ 
  h_{3,1}
  = &\ \theta_{1,0}\theta_{6,2}h_{2,0}+\theta_{1,1}(\theta_{6,1}
   + \theta_{6,5})h_{2,1}+\theta_{1,0}\theta_{6,4}h_{2,2} = h_{3,5},
\\
  h_{3,2}  
  = &\ \theta_{1,1}\theta_{6,4}h_{2,0}+\theta_{1,0}(\theta_{6,1}+\theta_{6,5})h_{2,1}+\theta_{1,1}\theta_{6,2}h_{2,2} = h_{3,4}, 
\\ 
  h_{3,3}  
  = &\ \theta_{1,0}\theta_{6,6}h_{2,0}+2\theta_{1,1}\theta_{6,3}h_{2,1}+\theta_{1,0}\theta_{6,0}h_{2,2}.\label{eq:h33}
\end{align}

From this expression, we can now give a formula for $C\Phi_6(q)$.
\begin{proposition}\label{prop:CPhi6}
The generating function $C\Phi_6(q)$ is equal to 
\begin{equation*} 
\frac{q^{\frac14}}{\eta(\tau)^6} \bigg(
   6\theta_{1,1}(\tau)^2\theta_{1,0}(\tau)\theta_{3,1}(3\tau)\theta_{2,1}(\tau) 
   + \theta_{1,0}(\tau)^3\big( \theta_{1,1}(6\tau)\theta_{1,1}(2\tau)+\theta_{1,0}(6\tau)\theta_{1,0}(2\tau) \big)
 \bigg) 
.
\end{equation*}
\end{proposition}

\begin{remark}
Note that Proposition~\ref{prop:CPhi6} agrees with the result of \cite{BS15}.  Indeed, using the identities
\begin{gather}
 \varphi(q) = \theta_{1,0}(\tau), \quad
 \qquad 2q^{\frac14}\psi\left(q^2\right) = \theta_{1,1}(\tau), \label{eq:BS} \\
 4q\psi(q)^3\psi\left(q^2\right)\psi\left(q^3\right) = \theta_{1,1}(\tau)\theta_{1,0}(\tau)\theta_{2,1}(\tau)\theta_{2,1}(3\tau), \label{eq:BS3}
\end{gather}
where 
\begin{equation}\label{eq:phipsi}
 \varphi(q):=\sum_{n\in\Z}q^{n^2} \quad \mbox{and} \quad \psi(q):=\sum_{n=0}^\infty q^{\frac{n(n+1)}{2}}
\end{equation}
is the notation used in \cite{BS15}, Proposition~\ref{prop:CPhi6} is equivalent to Theorem~2.1 of \cite{BS15}.  The identitities in \eqref{eq:BS} follow readily from the definitions, and \eqref{eq:BS3} can be checked using Sturm's Theorem.
\end{remark}

The following lemma is used in the proofs of Proposition~\ref{prop:CPhi6} and Proposition~\ref{prop:CPhi7}.  Using Sturm's Theorem, its proof is immediate, however, we include an elementary direct proof with the hope that it could be generalized.  Such generalizations may allow one to simplify or even give a closed form for $h_{k,k}$ for arbitrary $k\in \frac12 \Z_{\geq 0}$.

\begin{lemma}\label{lem:Theta3simplify}
We have 
\begin{align}
 \theta_{2,2}(\tau)\theta_{6,0}(\tau)+\theta_{2,0}(\tau)\theta_{6,6}(\tau)=&\ 2\theta_{2,1}(\tau)\theta_{6,3}(\tau), \label{eq:21lemA} \\
 \theta_{2,2}(\tau)\theta_{6,4}(\tau)+\theta_{2,0}(\tau)\theta_{6,2}(\tau)=&\ \theta_{2,1}(\tau)(\theta_{6,1}(\tau)+\theta_{6,5}(\tau)). \label{eq:21lemB}
\end{align}
\end{lemma}
 
\begin{proof}
Using the definitions and simplifying, we find that
 \[ \theta_{2,2}(\tau)\theta_{6,0}(\tau)+\theta_{2,0}(\tau)\theta_{6,6}(\tau)
 = \sum_{\substack{R,S\in \Z\\R\nequiv S \pmod{2}}} \(q^{\frac14}\)^{2R^2+6S^2}. \]

We now make the change of variables $s = \frac{R+S}{2},\ r = \frac{R-3S}{2}$.
Since $R$ and $S$ have opposite parity, $r,s\in \frac12+\Z$.  To ease the notation we also replace $q^{\frac14}$ with $q$.  Applying this change of variables gives 
\begin{equation*}
 \sum_{\substack{R,S\in \Z\\R\nequiv S \pmod{2}}} q^{2R^2+6S^2}
  = \sum_{\substack{r,s\in \frac12+\Z\\r\equiv s \pmod{2}}} q^{2r^2+ 6s^2}.
\end{equation*}  
We now write $r=1/2+j$ and $s=1/2+k$.  Since the condition $r\equiv s\pmod{2}$ is equivalent to $j\equiv k\pmod{2}$, it follows that this is equal to
\begin{align*}
 \sum_{\substack{j,k\in \Z\\j\equiv k \pmod{2}}} q^{2\(j+\frac12\)^2+6\(k+\frac12\)^2} 
  & = \sum_{j,k\in \Z} \(q^4\)^{2\(j+\frac14\)^2+6\(k+\frac14\)^2}
   + \sum_{j,k\in \Z} \(q^4\)^{2\(j-\frac14\)^2+6\(k-\frac14\)^2}.
\end{align*}
It is easy to see, replacing $q^4$ with $q$, that each of these final two sums is equal to $\theta_{2,1}(3\tau)\theta_{2,1}(\tau)$.  This proves \eqref{eq:21lemA}.  The identity \eqref{eq:21lemB} is similar; we omit the details.
\end{proof}

\begin{remark}
The terms 
\[ \theta_{2,0}\theta_{6,0}+\theta_{2,2}\theta_{6,2}\qquad\mbox{and}\qquad
    \theta_{2,0}\theta_{6,4}+\theta_{2,2}\theta_{6,2} \]
appear in the formulas for $C\Phi_6$ and $C\Phi_7$, and so it is natural to ask whether relations similar to \eqref{eq:21lemA} and \eqref{eq:21lemB} exist for these.  We do not, however, believe that such simple relations hold in these cases.
\end{remark}

\begin{proof}[Proof of Proposition~\ref{prop:CPhi6}]
From Theorem~2, we see that
\begin{align}\label{eq:CPhi6asCoeff}
 &C\Phi_6(q)=\mathrm{coeff}_{\left[\zeta^3\right]}\( \frac{\vartheta(z+\frac12;\tau)^6}{q^{\frac12}\eta(\tau)^6} \)
 =\frac{ \mathrm{coeff}_{\left[\zeta^3\right]}\vartheta_{3,3}(z;\tau)}{ q^{\frac12}\eta(\tau)^6 } h_{3,3}(\tau).
\end{align}
In order to obtain a more explicit formula for $C\Phi_6$, we plug the results of \eqref{eq:h2b} into \eqref{eq:h33}.  This yields
\begin{align}
 h_{3,3} &= 4\theta_{1,0}\theta_{1,1}^2\theta_{2,1}\theta_{6,3}
    + \theta_{1,0} \theta_{1,1}^2\big( \theta_{2,0}\theta_{6,6}+\theta_{2,2}\theta_{6,0}\big)
    + \theta_{1,0}^3 \big( \theta_{2,2}\theta_{6,6}+\theta_{2,0}\theta_{6,0}\big)  \nonumber
\\
 & = 6\theta_{1,0}\theta_{1,1}^2\theta_{2,1}\theta_{6,3}
    + \theta_{1,0}^3 \big( \theta_{2,2}\theta_{6,6}+\theta_{2,0}\theta_{6,0}\big), \label{eq:CPhi6Final}
\end{align}
using Lemma~\ref{lem:Theta3simplify} in the last step.  Since $\mathrm{coeff}_{[\zeta^3]}\vartheta_{3,3}(z;\tau)=q^{\frac34}$, the desired result follows from \eqref{eq:CPhi6asCoeff} and \eqref{eq:CPhi6Final}.
\end{proof}

\subsection{The case of $C\Phi_7$}

From Theorem~2 and \eqref{eq:thetamarelation}, after some simplification, we see that
\begin{align*}
 h_{\frac72,\frac72} = & \sum_{c\pmod{6}} h_{3,c} \theta_{21,7c-21} 
 = h_{3,0} \theta_{21,21} + 2h_{3,1}\theta_{21,14} + 2h_{3,2}\theta_{21,7} + h_{3,3}\theta_{21,0}.
\end{align*}
Using the expression for $h_{3,c}$ and $h_{2,c}$ given in equations~\eqref{eq:h2b}--\eqref{eq:h33}, we can write this completely in terms of $\theta_{p,q}$.  Upon simplifying, $h_{\frac72,\frac72}$ is equal to
\begin{align}
 &  \Big(\theta_{1,1}\theta_{6,0}\theta_{21,21} 
   + 2\theta_{1,0}\theta_{6,2}\theta_{21,14}
   + 2\theta_{1,1}\theta_{6,4}\theta_{21,7}
   + \theta_{1,0}\theta_{6,6}\theta_{21,0}\Big)
    \Big( \theta_{1,1}^2\theta_{2,0}+\theta_{1,0}^2\theta_{2,2} \Big) \nonumber \\
 & \quad +  
  4\theta_{1,0}\theta_{1,1}\theta_{2,1}
   \Big(\theta_{6,3}(\theta_{1,0}\theta_{21,21}+\theta_{1,1}\theta_{21,0})
   + (\theta_{6,1} + \theta_{6,5})
     (\theta_{1,1}\theta_{21,14}+\theta_{1,0}\theta_{21,7})\Big) \nonumber \\
 & \quad +  
  \Big(\theta_{1,1}\theta_{6,6}\theta_{21,21}
   + 2\theta_{1,0}\theta_{6,4}\theta_{21,14}
   + 2\theta_{1,1}\theta_{6,2}\theta_{21,7}
   + \theta_{1,0}\theta_{6,0}\theta_{21,0}\Big)
    \Big( \theta_{1,1}^2\theta_{2,2}+\theta_{1,0}^2\theta_{2,0} \Big). \label{eq:h7272}
\end{align}

With this in hand, we deduce the following result.

\begin{proposition}\label{prop:CPhi7}
We have $C\Phi_7(q)=h_{\frac72,\frac72}(\tau)/(q;q)_\infty^7$ with $h_{\frac72,\frac72}$ equal to
\footnotesize
 \begin{align*}
 &
6\theta_{1,0}\theta_{1,1}\theta_{2,1}
    \Big(\theta_{6,3}\big(\theta_{1,0}\theta_{21,21}+\theta_{1,1}\theta_{21,0}\big)
   + \big(\theta_{1,1}\theta_{21,14} + \theta_{1,0}\theta_{21,7}\big)
     \big(\theta_{6,1} + \theta_{6,5}\big)\Big) \\
 & \quad +
\Big(\theta_{1,0}^3\theta_{21,0}+\theta_{1,1}^3\theta_{21,21}\Big)
   \Big(\theta_{2,0}\theta_{6,0}+\theta_{2,2}\theta_{6,6} \Big) 
    + 2\Big( \theta_{1,0}^3\theta_{21,14} + \theta_{1,1}^3\theta_{21,7} \Big) 
   \Big(\theta_{2,0}\theta_{6,4}+\theta_{2,2}\theta_{6,2} \Big).
 \end{align*}
\normalsize
\end{proposition}
\begin{proof}
Note that $\mathrm{coeff}_{[\zeta^{7/2}]}\vartheta(z;\tau)=q^{7/8}$.  Thus, in complete analogy to \eqref{eq:CPhi6asCoeff}, we have 
 \[ C\Phi_7(q)= q^{\frac{7}{24}} \frac{h_{\frac72,\frac72}(\tau)}{\eta(\tau)^7} 
  = \frac{h_{\frac72,\frac72}(\tau)}{(q;q)_\infty^7}, \]
with $h_{\frac72,\frac72}$ as in \eqref{eq:h7272}.  In order to further simplify this expression, we first consider only the terms coming from the first and third lines of \eqref{eq:h7272}.  After rearranging and applying \eqref{eq:21lemA} and \eqref{eq:21lemB}, it can be shown that they are equal to
\footnotesize
\begin{align*}
& 2\theta_{1,0}\theta_{1,1}\theta_{2,1}\theta_{6,3}\Big( \theta_{1,0}\theta_{21,21}+\theta_{1,1}\theta_{21,0}\Big) 
   + \Big(\theta_{1,0}^3\theta_{21,0}+\theta_{1,1}^3\theta_{21,21}\Big)
   \Big(\theta_{2,0}\theta_{6,0}+\theta_{2,2}\theta_{6,6} \Big) \\
 & + 2\Big( \theta_{1,0}^3\theta_{21,14} + \theta_{1,1}^3\theta_{21,7} \Big) 
   \Big(\theta_{2,0}\theta_{6,4}+\theta_{2,2}\theta_{6,2} \Big) +
   2\theta_{1,0}\theta_{1,1}\theta_{2,1}\Big(\theta_{6,1}+\theta_{6,5} \Big)\Big(\theta_{1,0}\theta_{21,7}+\theta_{1,1}\theta_{21,14}\Big).
\end{align*}
\normalsize
Plugging these back into \eqref{eq:h7272} and simplifying gives the desired formula for $h_{\frac72,\frac72}$.
\end{proof}

At this point, one could give $C\Phi_7(q)$ in terms of $\eta$ and Klein forms using \eqref{eq:thetambasprod}.  Such a formula is quite long, however, so we refrain from writing it down explicitly.

\subsection{The case of $C\Psi_{4,0}$}
The formula for $C\Psi_{4,0}$ can be directly deduced by employing Theorem~2, yielding
\begin{align}\label{eq:CPsi40}
C\Psi_{4,0}(q)&=\frac{h_{2,0}}{q^{\frac{1}{2}}(q;q)_\infty^4} 
=\frac{h_{1,0}\theta_{2,0}h_{1,0}+h_{1,-1}\theta_{2,2}h_{1,1}}{q^{\frac{1}{2}}(q;q)_\infty^4} 
=\frac{\theta_{1,1}^2\theta_{2,0}+\theta_{1,0}^2\theta_{2,2}}{q^{\frac{1}{2}}(q;q)_\infty^4},
\end{align}
which will be utilized in Section~\ref{sec:congruences}.

\section{Congruences}\label{sec:congruences}

There are no doubt many congruences among $c\psi_{k,a}(n)$, and we have certainly not done an exhaustive search, but to illustrate some of the possibilities and to show how Theorem~2 can be used, we provide the following examples for $c\psi_{4,a}(n)$ with $a=0,1$.  Note that $a=2$ corresponds to $c\phi_4(n)$ which, as indicated in the introduction, has been extensively studied. We also note that the method of Jacobi forms this paper concerns allows one to identify that $C\Psi_{4,a}(q)$ are modular forms as well as their level and weight.  Thus a Sturm bound argument should suffice to prove the conjectural congruences. We have left this as open both because the naive Sturm bound one obtains is very large (though it can likely be substantially lowered by methods such as writing  in terms of vector-valued modular forms), and because this is merely an example of potential congruences.

\begin{conjthm}
The following congruences hold.

\begin{center}
\renewcommand{\arraystretch}{1.5}
\begin{tabular}{|c|c|}
\hline
Proved & Conjectured\\ \hline
$c\psi_{4,0}(2n+1)\equiv 0 \pmod{32}$ & $c\psi_{4,0}(7n+2)\equiv 0 \pmod{7}$\\
$c\psi_{4,0}(4n+3)\equiv 0 \pmod{64}$ & $c\psi_{4,1}(7n+3)\equiv 0 \pmod{7}$\\
$c\psi_{4,1}(n)\equiv 0 \pmod{4}$ & \\
\hline
\end{tabular}
\renewcommand{\arraystretch}{1}
\end{center}
\end{conjthm}

\begin{proof}[Proof that $c\psi_{4,0}(2n+1)\equiv 0 \pmod{32}$]
Note that $\phi(q)$ and $\psi(q)$ as given in \eqref{eq:phipsi} can be written as
\begin{equation*}
    \phi(q)=\sum\limits_{n=-\infty}^{\infty}q^{n^2}=(-q;q^2)_{\infty}^{2}(q^2;q^2)_{\infty} 
\end{equation*}
and
\begin{equation*}
    \psi(q)=\frac{1}{2}\sum\limits_{n=-\infty}^{\infty}q^{\frac{k(k-1)}{2}}=(-q;q^2)_{\infty}(q^4;q^4)_{\infty}.
\end{equation*}
Using the results of the previous section including \eqref{eq:BS}, the fact that
 \[ \theta_{2,0}(\tau) = \phi(q^2) \qquad \mbox{and} \qquad \theta_{2,2}(\tau) = q^{\frac{1}{2}}\psi(q^4), \]
and \eqref{eq:CPsi40}, we see that
\begin{align*}
      C\Psi_{4,0}(q) & = \frac{\theta_{1,1}^2\theta_{2,0}+\theta_{1,0}^2\theta_{2,2}}{q^{\frac{1}{2}}(q;q)_{\infty}^4}
      =
      \frac{4\psi^2(q^2)\phi(q^2)}{(q;q)_{\infty}^4}+\frac{2\phi^2(q)\psi(q^4)}{(q;q)_{\infty}^4}
      \\ & =
      \frac{4\psi^2(q^2)\phi(q^2)}{(q;q^2)_{\infty}^4(q^2;q^2)_{\infty}^4}+\frac{2\phi^2(q)\psi(q^4)}{(q;q^2)_{\infty}^4(q^2;q^2)_{\infty}^4}.
\end{align*}
Therefore,
\begin{align*}
    \sum\limits_{n=0}^{\infty}c\psi_{4,0}(2n+1)q^{2n+1} 
    &=
    \frac{1}{2}\left(\sum\limits_{n=0}^{\infty}c\psi_{4,0}(n)q^n-\sum\limits_{n=0}^{\infty}c\psi_{4,0}(n)(-1)^nq^n\right)
    \\ & =
    \frac{2\psi^2(q^2)\phi(q^2)}{(q^2;q^2)_{\infty}^4} \left(\frac{1}{(q;q^2)_{\infty}^4}-\frac{1}{(-q;q^2)_{\infty}^4}\right) 
    \\ & \qquad \qquad 
    + \frac{\psi(q^4)}{(q^2;q^2)_{\infty}^4}\left(\frac{\phi^2(q)}{(q;q^2)_{\infty}^4}-\frac{\phi^2(-q)}{(-q;q^2)_{\infty}^4}\right)
    \\ & =
    \frac{2\psi^2(q^2)\phi(q^2)}{(q^2;q^2)_{\infty}^4(q^2;q^4)_{\infty}^4}\left((-q;q^2)_{\infty}^4-(q;q^2)_{\infty}^4\right)
    \\ & \qquad \qquad 
    + \frac{\psi(q^4)}{(q^2;q^2)_{\infty}^2} \left(\frac{(-q;q^2)_{\infty}^4}{(q;q^2)_{\infty}^4}-\frac{(q;q^2)_{\infty}^4}{(-q;q^2)_{\infty}^4}\right).
\end{align*}
Note that in the final step we have used the fact that
\begin{equation}\label{phi}
\phi(q)=\frac{(q^2;q^2)_{\infty}^5}{(q;q)_{\infty}^2(q^4;q^4)_{\infty}^2}
\end{equation}
and
\begin{equation}\label{psi}
    \psi(q)=\frac{(q^2;q^2)_{\infty}^2}{(q;q)_{\infty}},
\end{equation}
which follow readily from \eqref{JTP}.

From Entry 25 of \cite[p.~40]{BR91}, we have
\begin{equation*}
    \phi^2(q)-\phi^2(-q)=8q\psi^2(q^4)
\end{equation*}
\begin{equation*}
    \phi^4(q)-\phi^4(-q)=16q\psi^4(q^2).
\end{equation*}
Applying these, we get
\begin{equation*}
    \sum\limits_{n=0}^{\infty}c\psi_{4,0}(2n+1)q^{2n+1}=
    \frac{16q\psi(q^4)\psi(q^2)}{(q^2;q^2)_{\infty}^{6}(q^2;q^4)_{\infty}^{4}}(\psi(q^2)\psi(q^4)\phi(q^2)+\psi^3(q^2)),
\end{equation*}
and by using \eqref{phi} and \eqref{psi} once again, we find that
\begin{equation*}
   \psi(q^2)\psi(q^4)\phi(q^2)+\psi^3{q^2}=2\frac{(q^4;q^4)_{\infty}^6}{(q^2;q^2)_{\infty}^3}.
 \end{equation*}
Plugging this back into the formula above, we see that
\begin{align}\label{cpsi4}
    \sum\limits_{n=0}^{\infty}c\psi_{4,0}(2n+1)q^{2n+1}&=32\frac{q\psi(q^4)\psi(q^2)}{(q^2;q^2)_{\infty}^{6}(q^2;q^4)_{\infty}^{4}}\frac{(q^4;q^4)_{\infty}^6}{(q^2;q^2)_{\infty}^3}
    \frac{32 q \psi(q^4)\psi^4(q^2)}{(q^2;q^2)_{\infty}^{6}(q^2;q^4)_{\infty}^{4}},
\end{align}
which immediately implies the desired result.
\end{proof}

In \eqref{cpsi4} we can divide both sides by $q$ and then replace $q^2$ by $q$.  Consequently, 
\begin{equation}\label{cpsi4qn}
    \sum\limits_{n = 0}^{\infty}c\psi_{4,0}(2n+1)q^n = \frac{32\psi(q^2)\psi^4(q^2)}{(q;q)_{\infty}^6(q^2;q^4)_{\infty}^{4}}
    = \frac{32(q^2;q^2)_{\infty}^7(q^4;q^4)_{\infty}^2}{(q;q)_{\infty}^{10}(q;q^2)_{\infty}^4},
\end{equation}
where in the final step we again used \eqref{psi}.

\begin{proof}[Proof that $c\psi_{4,0}(4n+3)\equiv 0 \pmod{64}$]
As a first step, we have
\begin{align*}
    \sum\limits_{n = 0}^{\infty}c\psi_{4,0}(4n+3)q^{2n+1}&=\frac{1}{2}\sum\limits_{n = 0}^{\infty}c\psi_{4,0}(2n+1)q^n-\sum\limits_{n = 0}^{\infty}c\psi_{4,0}(2n+1)(-1)^n q^n\\
    &=16\frac{(q^4;q^4)_{\infty}^2}{(q^2;q^2)_{\infty}^3}\left(\frac{1}{(q;q^2)_{\infty}^{14}}-\frac{1}{(-q;q^2)_{\infty}^{14}}\right)\\
    &=16\frac{(q^4;q^4)_{\infty}^2}{(q^2;q^2)_{\infty}^3(q^2;q^4)_{\infty}^{14}}(A^7-B^7),
\end{align*}
where $A=(-q;q^2)_{\infty}^2$, and $B=(q;q^2)_{\infty}^2$.  Since 
 \[ A^7-B^7 = (A-B)(A^6+A^5B+A^4B^2+A^3B^3+A^2B^4+AB^5+B^6)=:(A-B)f(A,B), \]
we find that 
\begin{align*}
    \sum\limits_{n=0}^{\infty}c\psi_{4,0}(4n+3)q^{2n+1}
    &=16\frac{(q^4;q^4)_{\infty}^2}{(q^2;q^2)_{\infty}^3(q^2;q^4)_{\infty}^{14}}(A-B)g(q)
\end{align*}
for a certain $q$-series $g(q)$ (where $g(q):=f(A,B)$).  Another application of Entry 25 of \cite[p.~40]{BR91} gives
\begin{equation}
    \phi(q)-\phi(-q)=4q\psi(q^8),
\end{equation}
hence, by \eqref{phi}, we have
\begin{equation}
    (A-B)(q^2;q^2)_{\infty} = (-q;q^2)_{\infty}^2(q^2;q^2)_{\infty}-(q;q^2)_{\infty}^2(q^2;q^2)_{\infty} = 4q\psi(q^8).
\end{equation}
Plugging this into the expression for $\sum c\psi_{4,0}(4n+3)q^{2n+1}$ above, the fact that $64|c\psi_{4,0}(4n+3)$ is evident.
\end{proof}

The proof that $c\psi_{4,1}(n)\equiv 0 \pmod{4}$ follows similar (but much easier) steps.  We leave it as an exercise for the interested reader.

\section{Proof of Theorem~3}\label{sec:motzkin}

\begin{figure}[h]\begin{center}
\begin{tikzpicture}
  \begin{axis}[
xtick={0,1,...,11},
ytick={0,1,...,15},
    separate axis lines,
    y axis line style= { draw opacity=0 },
    x axis line style= { draw opacity=0 },
    ymajorgrids,
		 xmajorgrids,
	   ]

    \addplot[name path=C,mark=none,thick,sharp plot] coordinates 
		{(0,0) (5,15) (6,15) (11,10)};
	
		\addplot[name path=D,mark=none,thick,sharp plot] coordinates
		{(3,9) (4,8) (5,9) (6,11) (7,12) (8,13)};
		
		\addplot[gray!50] fill between[of=C and D];
	



  \end{axis}
\end{tikzpicture}\end{center}
\caption{A $1$-shifted Motzkin path of length $11$ and rank $2$}\label{fig:MP-example2}
\end{figure}
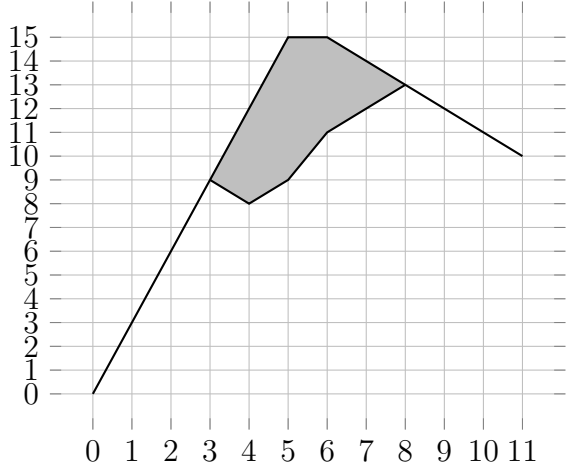

Recall that the definitions of $a$-shifted Motzkin paths of length $n$ and rank $\frac{k}{2}$ and the $(k,a)$-Motzkin path of length $n$ are provided in Definition~\ref{def:MotzkinPath}.  Given an $a$-shifted Motzkin path $w$, let $a_M(w)$ be the area under its graph.  Figure~\ref{fig:MP-example2} shows two $1$-shifted Motzkin paths of length $11$ and rank $2$.  The upper one $w_{\mathrm{max}}$ is that for which among all $(4,1)$-Motzkin paths $a_M(w_{\mathrm{max}})$ is maximal.  Notice that the area of the shaded region is
\begin{equation}\label{eq:bMw}
 b_M(w):=\frac{kn(n+1)}{2}+an(2n+1) - a_M(w),
\end{equation}
where $w$ is the lower path shown in Figure~\ref{fig:MP-example2}.

Recall that 
 \[ \widetilde{CM}_{k,a}^n(q):=\sum q^{a_M(w)}, \]
where the sum is over all $(k,a)$-Motzkin paths $w$ of length $n$, and thus
\begin{equation}\label{eq:CM_in_terms_of_bM}
 CM_{k,a}^{2n+1}(q):=q^{\frac{k}{2}n(n+1)+an(2n+1)}\widetilde{CM}_{k,a}^{2n+1}(q^{-1}) = \sum q^{b_M(w)}.
\end{equation}
We will prove below that for sufficiently large $n$ the set of paths of length $2n+1$ with fixed weight $b_M$ is in  bijection with the set of $(k,a)$-colored F-partition of weight $b_M$.  Moreover, in Lemma~\ref{lem:MPathFropPartBijection}, we will see that $b_M(w)$ is equal to the weight of the partition corresponding to $w$ under this bijection.  This implies, for bounded weights $b_M$ and $n$ sufficiently large, the coefficient of $q^{b_M}$ on the right hand side of \eqref{eq:CM_in_terms_of_bM} is precisely $c\psi_{k,a}(n)$.  In other words, letting $n\to\infty$,
 \[ \lim_{n\to \infty}CM_{k,a}^{2n+1}(q)=C\Psi_{k,a}(q), \]
which is the content of Theorem~3.  As outlined above, this requires us to establish the connection between Motzkin paths $w$ of length $2n+1$ and $(k,a)$-colored F-partitions.  Rather than work with Definition~\ref{def:kaFrobparts}, as discussed in the introduction, it will be easier to make this connection using the following alternative description of $(k,a)$-colored F-partitions.  

\begin{definition}\label{def:kaFrob-Drake}
Given $k\in \Z_{> 0}$ and $a\in \Z+\frac{k}{2}$ (with $a>0$), a \emph{$(k,a)$-generalized Frobenius partition} with weight $n:=r+s+\sum\limits_{\ell=0}^{r}a_\ell+\sum\limits_{j=0}^{s}b_j$ is a two-rowed array
\begin{equation} \label{eq:kaFrobPartition}
\begin{pmatrix}
a_1 & a_2 & \cdots  & a_r \\
b_1 & b_2 & \cdots  & b_s \\
\end{pmatrix}_T
\end{equation}
satisfying the following conditions.
\begin{itemize}
 \item[(i)] Each entry $a_j,b_j$ belongs to one of $k$ copies of non-negative integers.
 \item[(ii)] Each row is decreasing with respect to the lexicographical ordering.  (Meaning if $m_\ell=m\in\Z$ belongs to the $\ell$-th copy of $\Z$ and $n_j=n\in\Z$ belongs to the $j$-th copy of $\Z$, we have $m_\ell<n_j$ if and only if $m<n$ or $m=n$ and $\ell<j$.)
 \item[(iii)] The pair $(r,s)\neq(0,0)$ of nonnegative integers satisfies
\begin{equation}\label{eq:sconditions}
a-\frac{k}{2} \leq r-s \leq a+\frac{k}{2},
\end{equation}
 and $T$ is a subset of $\{1,2,\ldots,k\}$ of size $\# T = \frac{k}{2}-a+r-s$.
\end{itemize}
\end{definition}

\begin{lemma}
Definitions \eqref{def:kaFrobparts} and \eqref{def:kaFrob-Drake} are equivalent.
\end{lemma}
\begin{proof}
In order to see that the two definitions are equivalent, begin with a partition as in \eqref{eq:kaFrobPartition}.  Now add one to each $b_j$ for each $j=1,2,\cdots,s$, and extend the second row by $\# T$ of zeros, coloring each by the elements of $T$.  It is easy to see that this mapping gives a bijection between the two definitions.
\end{proof}

\begin{proof}[Proof of Theorem~3]
Given a $(k,a)$-colored generalized Frobenius partition as in Definition~\ref{def:kaFrob-Drake}:
\begin{equation}\label{eq:part}
\begin{pmatrix}
a_1 & a_2 & \cdots  & a_r \\
b_1 & b_2 & \cdots  & b_s \\
\end{pmatrix}_T,
\end{equation}
let $\alpha(p)$ be the number of times that $p$ appears in the first row of the partition, irrespective of coloring, i.e.,
 \[ \alpha(p) := \#\{ a_j\mid a_j = p, \ j=1,\ldots, r\}, \]
and let $U(p)$ be the set of all colorings of the elements $a_j$ with $a_j=p$.  Note that by part (ii) of Definition~\ref{def:kaFrob-Drake}, parts of the same color appear at most one, and thus $\# U(p) = \alpha(p)\leq k$.  We similarly define
 \[ \beta(q):= \#\{ b_j\mid b_j = q, \ j=1,\ldots, s\}, \]
and $V(q)$ is the set of colorings of all $b_j$ for which $b_j=q$.  By the same reasoning as above, $\#V(q) = \beta(q)\leq k$.

With these notations in place, for any $n\geq a_1,b_1$, we define a path of length $2n+1$ stipulated as follows:
\begin{itemize}
 \item[(i)] For each $p=0,1,\ldots,n-1$, the $(n-p)$-th step is $(1,\frac{k}{2}+a-\alpha(p))$ with coloring $\{1,2,\cdots,k\}\backslash U(p)$.
 \item[(ii)] The $(n+1)$-st step is $(1,r-s)$ with coloring $T$;
 \item[(iii)] For each $q=1,2,\ldots,n$, the $(n+1+q)$-th step is $(1,\beta(q)+a-\frac{k}{2})$ with coloring $V(q)$.
\end{itemize}
If $n>2\max\{a_1+1,b_1+1\}$, the path we have constructed does not go below the $x$-axis.  To see why, note that for each $p=n-1,n-2,\ldots,a_1+1$, it must be the case (since $a_1$ is the largest entry in the row) that $\alpha(p)=0$.  Therefore, every step $n-p=1,2,\ldots,n-a_1-1$ will be $(1,a+\frac{k}{2})$.  Beyond this point, the steepest decline that can occur would be if the following steps were each $(1,a-\frac{k}{2})$, which happens only if $\alpha(p)=k$ for all steps $n-p=n-a_1,n-a_1+1,\ldots,n$.  The $y$-coordinate of the path at the $n$-th step would then be
 \[ (n-a_1-1)\left(a+\frac{k}{2}\right) + (a_1+1)\left(a-\frac{k}{2}\right) = na + (n-2(a_1+1))\frac{k}{2}. \]
This is positive since $a\geq 0$ and $n>2(a_1+1)$.  A similar argument shows that for each of the steps $n,n+1,\ldots,2n+1$, the path remains above the $y$-axis.  Moreover, it is easily shown that the endpoint is $(2n+1,2na)$, hence the path defined above is a $(k,a)$-Motzkin path of length $2n+1$.

To a $(k,a)$-Motzkin path of length $2n+1$ we can append an additional step of size $(1,a+\frac{k}{2})$ at the start and one of size $(1,a-\frac{k}{2})$ at the end.  The resulting path, is also a $(k,a)$-Motzkin path, but now of length $2n+3$.  Using this, we define an equivalence relation on the set of $(k,a)$-Motzkin paths by setting two paths to be equivalent if one can be obtained from the other by successively performing the above procedure.  As recorded in Lemma~\ref{lem:MPathFropPartBijection} and proved below, every $(k,a)$-colored F-partition corresponds to a unique equivalence class of paths.  Additionally, since the quantity $b_M(w)$ given in \eqref{eq:bMw} is well-defined on equivalence classes, the discussion following \eqref{eq:CM_in_terms_of_bM} shows that Theorem~3 follows.
\end{proof}

\begin{lemma}\label{lem:MPathFropPartBijection}
The map on $(k,a)$-Frobenius partitions to equivalence classes of $(k,a)$-Motzkin paths is a bijection.  Moreover, under this mapping,
\begin{equation}\label{eq:bMwequality}
    b_M(w)=\frac{k}{2}n(n+1)+an(2n+1)-a_M(w)=r+s+\sum a_\ell+\sum b_j.
\end{equation}
That is to say that if $w$ is a $(k,a)$-Motzkin path of length $2n+1$ for some $n$ corresponding to a given partition, then $b_M(w)$ is equal to the weight of the partition.
\end{lemma}

\begin{figure}[h]\begin{center}
\begin{tikzpicture}
  \begin{axis}[
xtick={0,1,...,11},
ytick={0,1,...,15},
    separate axis lines,
    y axis line style= { draw opacity=0 },
    x axis line style= { draw opacity=0 },
    ymajorgrids,
		 xmajorgrids,
	   ]

    \addplot[name path=A,mark=none,thick,sharp plot] coordinates 
		{(0,0) (5,15) (6,15) (11,10)};
	
		\addplot[name path=B,mark=none,thick,sharp plot] coordinates
		{(3,9) (4,8) (5,9) (6,11) (7,12) (8,13)};
		
		\addplot[gray!50] fill between[of=A and B];
	
		\addplot[mark=none,sharp plot] coordinates
		{(5,9) (5,15)};

		\addplot[mark=none,sharp plot] coordinates
		{(6,11) (6,15)};

		\addplot[mark=none,sharp plot] coordinates
		{(4,8) (5,11)};

		\addplot[mark=none,sharp plot] coordinates
		{(6,13) (7,12)};
  \end{axis}
\end{tikzpicture}\end{center}
\caption{The $(4,1)$-Motzkin path of length 11 obtained from \eqref{eq:MPathPartition}}\label{fig:ThCproof}
\end{figure}
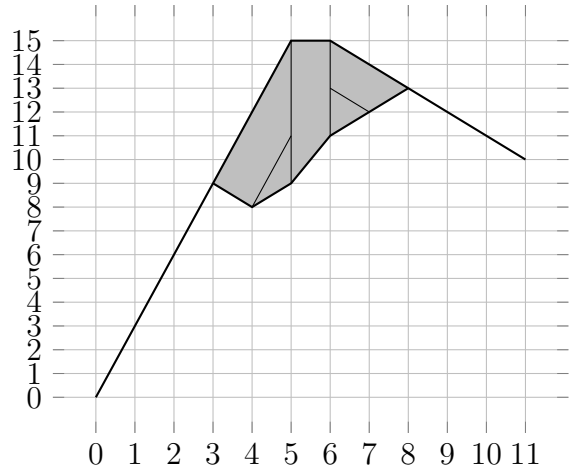

\begin{proof}
The fact that the map is a bijection is clear. That the equivalence relation does not effect the quantity $b_{M}$, is a consequence of the fact that the process of adding line segments of slopes $a\pm \frac{k}{2}$ on the left and right, respectively, merely shifts the figure, but has no effect on $b_M$, i.e., the area of the shaded region as in Figure~\ref{fig:ThCproof}.

To prove \eqref{eq:bMwequality}, let $(j,m_j)$ ($j=0,\ldots,2n+1)$) be the coordinates of the vertices of the Motzkin path $w$ corresponding to a partition \eqref{eq:part}.  For each interval $[j-1,j]$ with $j=1,2,\ldots,n$ we compute the area of the quadrilateral with vertices
 \[ \big(j-1,m_{j-1}\big), \quad \big(j,m_j\big), \quad \big(n,m_j+(n-j-1)(a+\ts{\frac{k}{2}})\big), \quad \big(n,m_{j-1}+(n-j)(a+\ts{\frac{k}{2}})\big). \]
If $j=n-p$, using property (i) above for the Motzkin path and shifting the lower left vertex to $(0,0)$, we find that the other three vertices are
 \[ \left(1,\ts{a+\frac{k}{2}} - \alpha(p)\right), \quad  \left(p + 1, (p + 1)\big(\ts{a+\frac{k}{2}}\big) - \alpha(p)\right), \left(p + 1, (p + 1)\big(\ts{a+\frac{k}{2}}-\alpha(p)\big)\right), \]
and the area of this quadrilateral is
 \[ p\alpha(p) + \frac{\alpha(p)}{2}. \]
Adding these all together gives
\begin{equation}
 \sum_{p=0}^{n-1} \Big( p\alpha(p) + \frac{\alpha(p)}{2}\Big) = \frac{r}{2} + \sum_{j=1}^r a_j.
\end{equation}
On the other hand, repeating this argument for the steps $[n+q,n+1+q]$ for each $q=1,2,\ldots,n$, we see that the corresponding areas add up to
\begin{equation}\label{eq:left}
 \sum_{q=1}^n \big(q\beta(q)+ \frac{\beta(q)}{2}\big) = \frac{s}{2} + \sum_{j=1}^s b_j.
\end{equation}\label{eq:right}
The remaining region over the interval $[n,n+1]$ is easily shown to have area $\frac{r+s}{2}$.  Hence, upon adding this to \eqref{eq:left} and \eqref{eq:right}, we find that the total area $b_M(w)$ is indeed equal to the weight of the partition, proving \eqref{eq:bMwequality}.
\end{proof}
 
To illustrate this lemma, we consider the case of $k=4$ and $a=1$, and the $(4,1)$ partition
\begin{equation}\label{eq:MPathPartition}
 \left(\begin{array}{cccccc} 1 & 1 & 1 & 1 & 0 & 0 \\ 1 & 1 & 0 & 0 \end{array}\right).
\end{equation}
Since $\alpha(p)$ and $\beta(q)$ are nonzero only in the cases $p,q=0,1$, we see that the slopes of each of the line segments between $0$ and $n-2$ are $a+\frac{k}{2}=3$.  Similarly, the slopes of the line segments starting at $n+3$ are each $a-\frac{k}{2}=-1$.  Finally, the slopes in between $n-2$ and $n+3$ are as shown in Figure~\ref{fig:ThCproof} (in which $n=5$).

If one computes the areas of each of the polygons in the shaded region of the figure, it is immediate that the area between $3$ and $5$ is indeed equal to $\frac{r}{2}+\sum a_j=5$; between $5$ and $6$ the area is $\frac{r+s}{2}=5$; between $6$ and $8$ the area is $\frac{s}{2}+\sum b_j = 4$.

\end{document}